\documentstyle [12pt,twoside, epsf]{article}

\makeatother

\def\picill#1by#2(#3)
{\vbox to #2
{\hrule width #1 height 0pt depth 0pt
\vfill\epsffile{#3}}}
 
\begin{document}

\date{}

\title{\Large\bf Remarks on Formal Knot Theory}
\author{Louis
H. Kauffman\\ Department of Mathematics, Statistics \\ and Computer Science (m/c
249)    \\ 851 South Morgan Street   \\ University of Illinois at Chicago\\
Chicago, Illinois 60607-7045\\ $<$kauffman@uic.edu$>$\\}

\maketitle

\thispagestyle{empty}

\section{Introduction}
This article is a supplement to the text of Formal Knot Theory.
It can be used as a companion and guide to the original book, as well as an 
introduction to some of the subsequent developments in the the theory of knots. 
Here is a table of contents for these Remarks.
\bigbreak

\begin{enumerate}
\item Alexander polynomial
\item Clock Theorem
\item The FKT state sum and its properties
\item Coloring edges and regions
\item Remark on the Duality Conjecture
\item Bracket Polynomial and Jones polynomial
\item Remarks on quantum topology and Khovanov Homology
\end{enumerate}

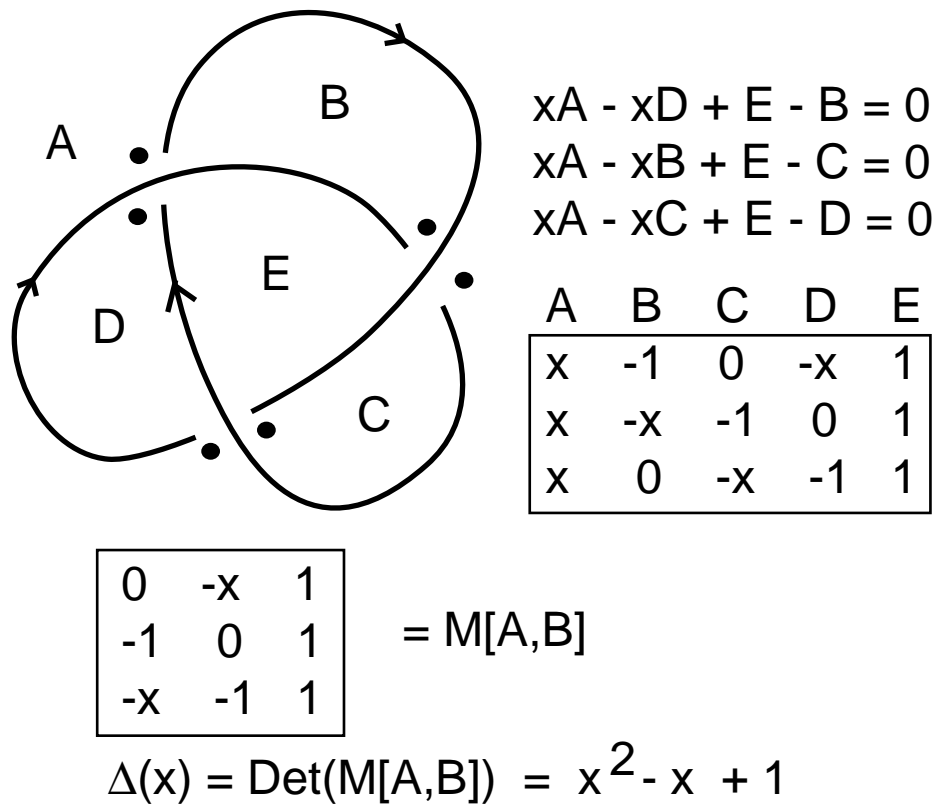
\section{The Alexander Polynomial}
In this section we give a short description of the Alexander Polynomial as it was defined by J. W. Alexander in
his 1928 paper \cite{Alex}. Alexander defined the polynomial as the determinant of a matrix associated with 
an oriented link diagram. In his paper he shows that his
determinant is well-defined and invariant under Reidemeister moves up to sign and powers of the variable $x.$ This method of presentation makes the Alexander
polynomial accessible to anyone who knows basic linear algebra and is willing to learn the basics of knot and link diagrams. 
\bigbreak

Later in the paper, Alexander describes how the matrix whose determinant yields the polynomial is related to the fundamental group of the complement of 
the knot or link. This relationship is described in slightly more modern language in Formal Knot Theory Appendix 1, and we shall not discuss it in this
introduction. Suffice it to say that the Alexander matrix is a presentation of the abelianization of the commutator subgroup of the fundamental group as 
a module over the group ring of the integers. Since Alexander's time, this relationship with the fundamental group has been understood very well.
For more information about this point of view we recommend that the reader consult \cite{Fox1,Fox2,Rolfsen,OnKnots}. See also the
appendix to Formal Knot Theory for a derivation of the structure of this module that corresponds to the algorithm in Alexander's original paper.
\bigbreak

Alexander gave a useful notation for generating the matrix associated to a diagram. At each crossing two dots are placed just to the left of the
undercrossing arc, one before and one after the overcrossing arc at the crossing. Here one views the crossing so that the undercrossing arc is vertical
and the overcrossing arc is horizontal. See Figure 1. 
\bigbreak

Four regions meet locally at a given crossing. Letting these be labeled generically $\{ A,B,C,D \},$ as shown in Figure 1, Alexander associates the 
equation $$xA - xB + C - D = 0$$ to that crossing. Here $A,B,C,D$ proceed cyclically around the crossing, starting at the top dot. In this way the two
regions containing the dots give rise to the two occurrences of $x$ in the equation. If some of the regions are the same at the crossing, then the equation is
simplified by that equality. For example, if $A=D$ then the equation becomes $xA - xB + C -A = 0.$ Each node in a diagram $K$ gives an equation 
involving the regions of the diagram. Alexander associates a matrix $M_{K}$ whose rows correspond to the nodes of the diagram, and whose columns correspond
to the regions of the diagram. Each nodal equation gives rise to one row of the matrix where the entry for a given column is the coefficient of that 
column (understood as designating a region in the diagram) in the given equation. If $R$ and $R'$ are adjacent regions, let $M_{K}[R,R']$ denote the
matrix obtained by deleting the corresponding columns of $M_{K}.$ Finally, define the {\em Alexander polynomial} $\Delta_{K}(x)$ by the formula
$$\Delta_{K}(x) \dot{=} Det(M_{K}[R,R']).$$

$$ \picill5inby3in(F1)  $$
\begin{center}
{\bf Figure 1  - Alexander Labeling. }
\end{center}

\noindent The notation $ A \dot{=} B$ means that $A = \pm x^{n}B$ for some integer $n.$ Alexander proves that his polynomial is well-defined, independent of
the choice of adjacent regions and invariant under the Reidemeister moves  up to $\dot{=}.$ The proof in Alexander's paper is elementary and combinatorial.
We shall not repeat it here.
\bigbreak

In Figure 2 we show a the calculation of the Alexander polynomial of the trefoil knot using this method.

$$ \picill5inby5in(F2)  $$
\begin{center}
{\bf Figure 2  - The Alexander Polynomial. }
\end{center}

\noindent In this figure we show the diagram of the knot, the labelings and the resulting full matrix and the square matrix resulting from deleting
two columns corresponding to a choice of adjacent regions. Computing the determinant, we find that the the Alexander polynomial of the trefoil knot is given
by the equation $\Delta \dot{=} x^2 - x + 1.$ 

\section{Reformulating the Alexander Polynomial as a State Summation}
Formal Knot Theory is primarily about a reformulation of the Alexander polynomial as a state summation. This means that we shall give a formula for 
the Alexander polynomial that is a sum of evaluations of combinatorial configurations related to the knot or link diagram. The {\em states} are these 
combinatorial configurations. These states are directly related to the expansion of the determinant that defines the 
polynomial. Graph-theoretic miracles occur, and it turns out that one can get a fully normalized version of the 
Conway version of the Alexander polynomial by using these states. In this section we shall show how the states emerge naturally from 
Alexander's determinant.
\bigbreak

Given a square $n \times n$ matrix $M_{ij}$, we consider the expansion formula for the determinant of $M:$
$$Det(M) = \Sigma_{\sigma \in S_{n}}sgn(\sigma)M_{1\sigma_{1}}\cdots M_{n\sigma_{n}}.$$
Here the sum runs over all permutations of the indices $\{ 1,2, \ldots ,n \}$ and $sgn(\sigma)$ denotes the sign of a given permutation $\sigma.$
In terms of the matrix, each product corresponds to a {\em choice} by each column of a single row such that each row is chosen exactly once. The order of 
rows chosen by the columns (taken in standard order) gives the permutation whose sign is calculated.
\bigbreak

Consider our description of Alexander's determinant as given in the previous section. Each crossing is labeled with Alexander's dots so that we know that 
the four local quadrants at a crossing are each labeled with $x$, $-x$, $1$ or $-1.$ The matrix has one row for each crossing and one column for each region.
Two columns corresponding to adjacent regions $A$ and $B$ are deleted from the full matrix to form the matrix $M[A,B]$, and we have the Alexander polynomial
$\Delta_{K}(X) \dot{=} Det(M[A,B]).$ 
\bigbreak

In the Alexander determinant expansion the {\em choice} of a row by a column corresponds to {\em a region choosing a node} in the link diagram.
The only nodes that a region can choose giving a non-zero term in the determinant are the nodes in the boundary of the given region. Thus the terms in the
expansion of $Det(M[A,B])$ are in one-to-one correspondence with decorations of the flattened link diagram (i.e. we ignore the over and under crossing
structure) where each region (other than the two deleted regions corresponding to the two deleted columns in the matrix) labels one of its nodes. We call
these  labeled flat diagrams the {\em states} of the original link diagram. See Figure 3 for a list of the states of the trefoil knot.
In this figure we show the states and the corresponding matrix forms with columns choosing rows that correspond to each state.
\bigbreak

$$ \picill5inby3in(F3)  $$
\begin{center}
{\bf Figure 3  - States with Markers }
\end{center}

$$ \picill2inby1.5in(F4)  $$
\begin{center}
{\bf Figure 4  - A Black Hole }
\end{center}

$$ \picill5inby6in(F5)  $$
\begin{center}
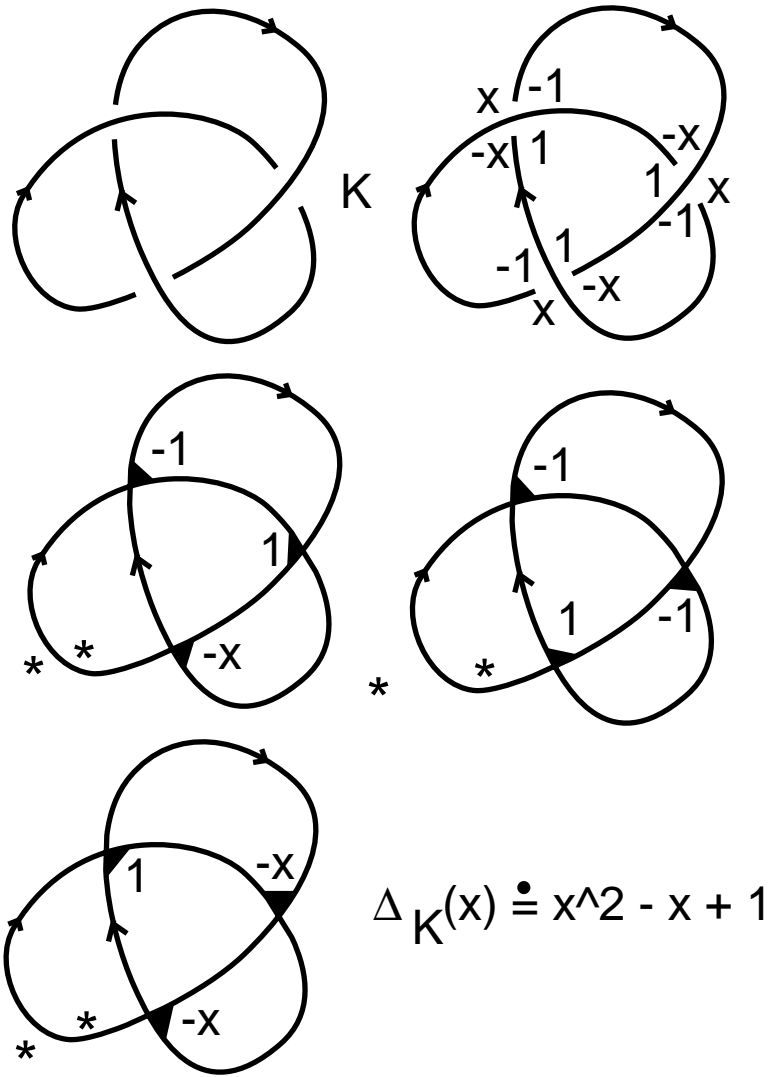
{\bf Figure 5  - State Sum Calculation of Alexander Polynomial }
\end{center}

It is useful to have terminology for a flattened link diagram. We call such a flattened
diagram a {\em link universe} or just {\em universe} for short. Such a diagram is a graph in the plane with four edges
incident to each vertex. The vertex carries no information about under or over crossings of curves as in a knot or link diagram.  Thus a universe is a $4$-regular planar graph. We say that a universe is {\em connected} if
it is connected as a planar graph. On the other hand, we say that a link diagram has $k$ components if it represents an embedding of $k$ circles in three
dimensional space. One counts the number of components of a link diagram by walking along the diagram and crossing at each vertex, counting the number of
cycles needed to use all the edges in the diagram. In the corresponding flattened diagram the operation of crossing at a vertex means that, at a vertex, one
chooses to continue the walk along the unique edge that is {\em not} adjacent to the edge one is traversing. The planar embedding of the graph defines this
adjacency. Thus we can speak of the number of {\em link components of a universe}. This number can be greater than one even when the universe is
connected. Two circles intersecting transversely in two points form a connected universe with two link components. Note that given a universe $U$ with $n$
nodes there are $2^{n}$ possible link diagrams that can be made from $U$ by choosing a crossing at each node.
\bigbreak

At this point we have almost a full combinatorial description of Alexander's determinant. The only thing missing is the permutation signs.
One can pick up the permutations from the state labeling, but there is a better way. Call a state marker (label at a node as shown in Figure 3)
a {\em black hole} if it labels a quadrant where both oriented segments point toward the node. See Figure 4 for an illustration of this concept.
\bigbreak

\noindent Let $S$ be a state of the diagram K. Consider the parity $$(-1)^{b(S)}$$ where $b(S)$ is the number of black holes in
the state
$S$. Then it turns out that  up to one global sign $\epsilon$ depending on the ordering of nodes and regions, we have $$(-1)^{b(S)} = \epsilon \,\,
sgn(\sigma(S))$$ where 
$\sigma(S)$ is the permutation of nodes induced by the choice of ordering of the regions of the state. This gives a purely diagrammatic access to 
the sign of a state and allows us to write
$$\Delta_{K}(x) \dot{=} \Sigma_{S}<K|S>(-1)^{b(S}$$ where $S$ runs over all states of the diagram for a given choice of deleted
adjacent regions, and $<K|S>$ denotes the product of the Alexander nodal labels at the quadrants indicated by the state labels in the state $S.$ 
We call $<K|S>$ the {\em product of the vertex weights}.
Thus we have a precise reformulation of the Alexander polynomial as a state summation.
\bigbreak

In Figure 5 we illustrate the calculation of the Alexander polynomial of the trefoil knot using this state summation. Here we show the contributions of
each state to a product of terms and in the polynomial we have followed the state summation by taking into account the number of black holes in each state. 
The most mysterious thing about
this state sum is the harmony of signs $$(-1)^{b(S)} = \epsilon \,\, sgn(\sigma(S)).$$ We shall explain this harmony in the next section via the {\em
Clock Theorem}. In Section 5 we shall reformulate the Alexander nodal weights to obtain the Conway polynomial.
\bigbreak

\section{The Clock Theorem}
In this section we continue the discussion of the states of a knot or link universe that we began in the previous section.
Consider the relationship of the two states shown in Figure 6. As the Figure shows, one state is obtained from the other by exchanging the markers at
two nodes so that in the state $S$ we have a marker at $i$ in region $A$ and a marker at $j$ in region $B$, while in state $S'$ we have a marker at
$i$ in region $B$ and a marker at $j$ in region $A.$ The prerequisite for this kind of exchange is that the regions $A$ and $B$ are adjacent with the nodes
$i$ and $j$ incident to both regions. Note also that one can think of this exchange as mediated by rotating each marker in the same clock direction by
ninety degrees. Thus $S'$ is obtained from $S$ by one clockwise rotation. We call such an exchange a {\em clocking move}. The Clock Theorem
(See \cite{FKT}) states that {\em any two states of a connected link universe are connected by clocking moves.} 
\bigbreak

$$ \picill3inby3in(F6)  $$
\begin{center}
{\bf Figure 6  - A Clockwise Move }
\end{center}

It is easy to see 
that the parity of the number of black holes in a state (see the previous section for the definition of a black hole) is {\em changed} by any single clocking
move. This means that if $S'$ is obtained from $S$ by one clocking move, then $sgn(S') = - sgn(S)$ where $sgn(S) = (-1)^{b(S)}$ as in the previous
section. On the other hand, if we consider the permutation of nodes $\sigma(S)$ that is associated to a given state $S$ then it is clear that 
$\sigma(S')$ and $\sigma(S)$ differ by one transposition and therefore $sgn(\sigma(S')) = - sgn(\sigma(S)).$ This verifies the predetermined 
harmony between the permutation signs coming from Alexander's determinant and the state signs coming from the black holes. Thus the state summation
model for the Alexander polynomial described in the previous section is seen to be founded on the Clock Theorem's statement that any two states are
connected by a series of clocking moves.
\bigbreak

An important combinatorial structure related to a link diagram is the {\it checkerboard graph}. This is a graph derived from a {\it checkerboard shading}
of the diagram. A checkerboard shading is a coloring of the regions of the diagram in two colors so that two adjacent regions receive different colors.
This can always be accomplished. Try this as an exercise. 
\bigbreak

$$ \picill4inby6.5in(F7)  $$
\begin{center}
{\bf Figure 7  - The Checkerboard Graph. }
\end{center}

\noindent See Figure 7 for an illustration of a checkerboard coloring. 
We usually refer to the two types of colored region as {\it shaded} and {\it unshaded}.The checkerboard graph of the diagram
is a graph with one node for each shaded region and one edge for each crossing between two shaded regions.
The states admit useful reformulations. They are in one-to-one correspondence with the maximal trees in the checkerboard graph of the knot or link.
They are also in one-to-one correspondence with Jordan-Euler trails on the link diagram. For the definitions of the checkerboard graph and the concept of 
Jordan-Euler trail, see pages 13 to 20 of Formal Knot Theory. After reading this, the
reader may enjoy the following exercise:

\noindent {\bf Exercise.} Let $K$ be a connected alternating link diagram. Prove that the absolute value of $\Delta_{K}(-1)$ is equal to the number of maximal
trees  in the checkerboard graph of $K$. (Hint: Show that when $x = -1$ every state contributes the same sign to the state summation for $\Delta_{K}(-1).$)
\bigbreak

\noindent{\bf Remark.} Ten years after the publication of Formal Knot Theory, a planar graph version of the Clock Theorem was obtained independently by
James Propp \cite{Propp}. It remains to be seen how his methods reflect on the knot theory that is associated with this result.
\bigbreak

\section{Reformulating the State Sum}
Recall that four regions meet locally at a given crossing. Letting these be labeled generically $\{ A,B,C,D \},$ as shown in Figure 1,
Alexander associates the  equation $xA - xB + C - D = 0$ to that crossing. From this we obtained the vertex weights for our state sum for the 
Alexander polynomial. We now point out that if we change the weights to correspond to the equation $xA + xB + C + D = 0$ (that is we remove all the 
negative signs from Alexander's labels), the the polynomial $\Delta_{K}(x)$ will only change by a global sign. To see this examine the effect of the signs
in Alexander's weights on a single clock move. It is easy to see that the product of these signs does not change when two states are related by a clock
move. Thus the signs in the Alexander weights only contribute a global sign to the polynomial. We therefore remove them. See Figure 8.
\bigbreak

$$ \picill3inby2in(F8)  $$
\begin{center}
{\bf Figure 8  - New Alexander Labeling. }
\end{center}

$$ \picill5inby4in(F9)  $$
\begin{center}
{\bf Figure 9  - Changing the Alexander Labeling. }
\end{center}

Having removed the signs from the Alexander weights, we show the new weights in Figure 8. In Figure 9
we make another change by replacing $x$ by $x^2$ and then writing the result on the weights that corresponds to multiplying each local set of 
weights by $1/x.$  Using the $x^2$ substitution, the
polynomial has changed by a power of
$x.$ Finally, the weights $x$ and $1/x$ that occur between  oppositely oriented arcs are removed. This again changes the polynomial by a factor that is a
power of
$x$ (checked again by using the Clock Theorem). We end up with the weights in Figure 10 after exchanging $x$ and $1/x$. These weights give us a state sum
$\Omega_{K}(x)$ with
$$\Omega_{K}(\sqrt{x}) \dot{=}
\Delta_{K}(x).$$
\bigbreak

$$ \picill5inby2.5in(F10)  $$
\begin{center}
{\bf Figure 10  - Conway Polynomial Weights. }
\end{center}

In Figure 10 the new weights (shown in a box) are called the Conway polynomial weights because the state sum now computes the Conway normalized version of
the  Alexander polynomial. The Conway normalization no longer has any change of sign or power of the variable under the Reidemeister moves, and it satisfies 
a ``skein relation" as described below.
\bigbreak

$$ \picill5inby2in(F11)  $$
\begin{center}
{\bf Figure 11  - The Skein Triple. }
\end{center}

Let three oriented link diagrams $K_{+}, K_{-}, K_{0}$ be related by changes at the site of a single crossing so that 
$K_{+}$ has a positive crossing, the crossing is switched in $K_{-}$ and smoothed in $K_{0}$ as in Figure 11. Then in Formal Knot Theory we show
that $$\Omega_{K_{+}} - \Omega_{K_{+}} = (x - x^{-1})\Omega_{K_{0}}.$$
It follows from this identity that $\Omega_{K}$ is a polynomial in $z = x - x^{-1}.$ and so we usually write $\Omega_{K}(z)$ rather than $\Omega_{K}(x).$
It is also that case that $\Omega_{K}$ is strictly invariant under all the Reidemeister moves: If $K$ and $K'$ are equivalent link diagrams, that
$\Omega_{K}(z) = \Omega_{K'}(z).$ There is no need to use the $equals \, dot$ relation any longer. The state sum for $\Omega_{K}(z)$ gives a combinatorial
model for the Conway version \cite{Conway} of the Alexander polynomial. 
\bigbreak

There are other ways to model the Conway-Alexander Polynomial. A method using an orientable spanning surface for the link and a matrix of linking numbers of
curves on this surface (the Seifert matrix) is discussed in \cite{OnKnots}.
\bigbreak

$$ \picill5inby6in(F12)  $$
\begin{center}
{\bf Figure 12  - Conway Calculation. }
\end{center}

In Figure 12 we illustrate the calculation of this state sum for the trefoil knot. The state summation model for the Conway - Alexander polynomial is one of
the main points of the book Formal Knot Theory. We use this model to prove results about alternating and alternative knots. The state summation gives
access to properties of the Alexander polynomial that are difficult to see from its definition as a determinant. This state 
summation model paved the way to discovering the bracket state sum model for the Jones polynomial \cite{Jones,Bracket}. The bracket model was discovered after
Formal Knot Theory had already been published in 1983. In this new edition, we include a paper \cite{New} on the bracket polynomial so that the reader can see
the  continuity between these two state summation models.
\bigbreak

\section{Coloring Edges and Regions}
A well-known \cite{OnKnots,Fox2} way to show that some knots are knotted is the method of Fox coloring. In this method each arc of the knot diagram
is labeled with an element of $Z/NZ$ for an appropriate modulus $N.$ At each crossing we require the equation $$a + c = 2b$$ where $a,b,c$ are the labels
(colors) from $Z/NZ$ incident to the crossing, and $b$ is the label of the overcrossing arc, while $a$ and $c$ are the labels of the undercrossing arcs that
meet the crossing. If a knot can be colored non-trivially, then it is not hard to see that it is knotted by examining how colorings 
change under the Reidemeister moves. In fact one sees that each Reidemeister move induces a unique change in the coloring of a knot or link, and that for
knots, a coloring with more than two colors is preserved (in the sense of continuing to have more than two colors) under Reidemeister moves that make only
local changes in the diagram and in the coloring. 
\bigbreak

$$ \picill5inby6in(F13)  $$
\begin{center}
{\bf Figure 13  - Coloring Faces and Edges }
\end{center}

It is the purpose of this section to relate the edge-coloring method to a face-coloring method that is very close to the structure of the FKT model of the 
Alexander polynomial. To see this relationship, view Figure 13. Here we see a labeling of faces incident to a crossing. Suppose that that these faces are
labeled from $Z/NZ$ and that they satisfy the equation $$A+B = C +D$$ where pairs $A,B$ and $C,D$ labels for faces sharing the over-crossing line.
Now associate to each arc at the crossing, the sum of the labels on either side of it. Thus
$$a = B + C,$$
$$b= A + B = C + D,$$
$$c = A + D.$$
Then we see that 
$$a + c = B + C + A + D = A + B + C + D = 2b.$$
Thus a labeling of faces satisfying the $A + B = C + D$ rule gives rise to a Fox coloring of the arcs of the diagram. One can go back and forth this
way between arc colorings and face colorings. 
\bigbreak

The next observation is to realize that the rule $A + B = C + D$ is none other than the equation at the crossing that corresponds to the 
Alexander polynomial with $x = -1$. To see this, use the revised Alexander labeling of Figure 8 and the discussion related to this figure.
It is then not hard to see that if one takes as modulus $N = \Delta_{K}(-1)$, then there exist colorings of the faces of the diagram that satisfy
these equations. The absolute value of $\Delta_{K}(-1)$ is called the {\it determinant} of the knot $K,$ and is often denoted
$$Det(K) = Abs(\Delta_{K}(-1)).$$
The determinant of the knot is itself an invariant of the knot, and any Fox coloring modulus must divide it.
\bigbreak

We have given an exercise (Section 4) to show that $Det(K)$ is equal to the number of maximal trees in the checkerboard graph of the knot $K$ when 
$K$ is alternating. The means that the maximal tree count gives the color modulus in this case.
\bigbreak

\section{The Duality Conjecture}
This section is a note to inform the reader that the Duality Conjecture stated on page $57$ of Formal Knot Theory has been proved by 
Gilmer and Litherland \cite{GL}. Of course with that encouragement, the reader may enjoy finding a proof for herself without consulting the literature.
\bigbreak

\section{The Bracket Polynomial and The Jones Polynomial}
The bracket model of the Jones polynomial \cite{Bracket} is a state summation model whose structure is very close to the FKT model of the
Alexander polynomial. Since a paper about the bracket model \cite{New} is included with this reprinting of Formal Knot Theory, we will only briefly sketch the
bracket model here. One point worth making is that one can re-write the bracket model as a summation over all the trees in the checkerboard graph \cite{T,New}
or, equivalently,  over all the  single cycle states of the diagram. This means that the very same set of states that yields the Alexander-Conway polynomial
can be used to produce the  Jones polynomial. A deeper understanding along these lines of the relationship of the Clock Theorem to the Jones polynomial is an
open question.
\bigbreak

It is an open problem whether there exist classical knots (single component loops) that are knotted and yet have unit Jones polynomial.
In other words, it is an open problem whether the Jones polynomial can detect all knots. 
There do exist families of links whose linkedness is undetectable by the Jones polynomial \cite{Morwen,EKT}. 
\bigbreak

The {\em bracket polynomial} , $<K> \, = \, <K>(A)$,  assigns to each unoriented link diagram $K$ a 
Laurent polynomial in the variable $A$, such that
   
\begin{enumerate}
\item If $K$ and $K'$ are regularly isotopic diagrams, then  $<K> \, = \, <K'>$.
  
\item If  $K \amalg O$  denotes the disjoint union of $K$ with an extra unknotted and unlinked 
component $O$ (also called `loop' or `simple closed curve' or `Jordan curve'), then 

$$< K \amalg O> \, = \delta<K>,$$ 
where  $$\delta = -A^{2} - A^{-2}.$$

\item $<K>$ satisfies the following formulas 

$$<\mbox{\large $\chi$}> \, = A <\mbox{\large $\asymp$}> + A^{-1} <)(>$$
$$<\overline{\mbox{\large $\chi$}}> \, = A^{-1} <\mbox{\large $\asymp$}> + A <)(>,$$
\end{enumerate}

\noindent where the small diagrams represent parts of larger diagrams that are identical except  at
the site indicated in the bracket. We take the convention that the letter chi, \mbox{\large $\chi$},
denotes a crossing where {\em the curved line is crossing over the straight
segment}. The barred letter denotes the switch of this crossing, where {\em the curved
line is undercrossing the straight segment}. 
\bigbreak

In computing the bracket, one finds the following behaviour under Reidemeister move I: 
  $$<\mbox{\large $\gamma$}> = -A^{3}<\smile> \hspace {.5in}$$ and 
  $$<\overline{\mbox{\large $\gamma$}}> = -A^{-3}<\smile> \hspace {.5in}$$

\noindent where \mbox{\large $\gamma$}  denotes a curl of positive type as indicated in Figure 14, 
and  $\overline{\mbox{\large $\gamma$}}$ indicates a curl of negative type, as also seen in this
figure. The type of a curl is the sign of the crossing when we orient it locally. Our convention of
signs is also given in Figure 14. Note that the type of a curl  does not depend on the orientation
we choose.  The small arcs on the right hand side of these formulas indicate
the removal of the curl from the corresponding diagram.  

\bigbreak
  
\noindent The bracket is invariant under regular isotopy and can be  normalized to an invariant of
ambient isotopy by the definition  
$$f_{K}(A) = (-A^{3})^{-w(K)}<K>(A),$$ where we chose an orientation for $K$, and where $w(K)$ is 
the sum of the crossing signs  of the oriented link $K$. $w(K)$ is called the {\em writhe} of $K$. 
The convention for crossing signs is shown in  Figure 14. The original Jones polynomial, $V_{K}(t)$ \cite{Jones} is then given by the formula
$$V_{K}(t) = f_{K}(t^{-1/4}).$$

$$ \picill4.5inby3in(F14) $$
\begin{center} {\bf Figure 14 - Crossing Signs and Curls} 
\end{center}
\vspace{3mm}

\noindent One useful consequence of these formulas is the following {\em switching formula}
$$A\mbox{\large $\chi$} - A^{-1} \overline{\mbox{\large $\chi$}} = (A^{2} - A^{-2})\mbox{\large $\asymp$}.$$ Note that 
in these conventions the $A$-smoothing of $\mbox{\large $\chi$}$ is $\mbox{\large $\asymp$},$ while the $A$-smoothing of
$\overline{\mbox{\large $\chi$}}$ is $><.$ Properly interpreted, the switching formula above says that you can switch a crossing and 
smooth it either way and obtain a three diagram relation. This is useful since some computations will simplify quite quickly with the 
proper choices of switching and smoothing. Remember that it is necessary to keep track of the diagrams up to regular isotopy (the 
equivalence relation generated by the second and third Reidemeister moves). Here is an example. View Figure 15.
\bigbreak

\noindent You see in Figure 15, a trefoil diagram $K$, an unknot diagram $U$ and another unknot diagram $U'.$  Applying the switching formula,
we have $$A^{-1} <K> - A <U> = (A^{-2} - A^{2}) <U'>$$ and 
$<U> = -A^{3}$ and $<U'> =(-A^{-3})^2 = A^{-6}.$  Thus $$A^{-1} <K> - A(-A^{3}) = (A^{-2} - A^{2}) A^{-6}.$$ Hence
$$A^{-1} <K> = -A^4 + A^{-8} - A^{-4}.$$  Thus $$<K> = -A^{5} - A^{-3} + A^{-7}.$$ This is the bracket polynomial of the trefoil diagram $K.$
\bigbreak

\noindent Since the trefoil diagram $K$ has writhe $w(K) = 3,$ we have the normalized polynomial 
$$f_{K}(A) = (-A^{3})^{-3}<K> = -A^{-9}(-A^{5} - A^{-3} + A^{-7}) = A^{-4} + A^{-12} - A^{-16}.$$ The asymmetry of this polynomial under the 
interchange of $A$ and $A^{-1}$ proves that the trefoil knot is not ambient isotopic to its mirror image.
\bigbreak

$$ \picill3inby2in(F15) $$
\begin{center} {\bf Figure 15 -- Trefoil and Two Relatives} \end{center}
\bigbreak

The bracket model for the Jones polynomial is quite useful both theoretically and in terms
 of practical computations. One of the neatest applications is to simply compute, as we have done, $f_{K}(A)$ for the
trefoil knot $K$ and determine that  $f_{K}(A)$ is not equal to $f_{K}(A^{-1}) = f_{-K}(A).$  This
shows that the trefoil is not ambient isotopic to its mirror image, a fact that is much harder to
prove by classical methods.
\bigbreak

\noindent {\bf The State Summation.} In order to obtain a closed formula for the bracket, we now describe it as a state summation.
Let $K$ be any unoriented link diagram. Define a {\em state}, $S$, of $K$  to be a choice of
smoothing for each  crossing of $K.$ There are two choices for smoothing a given  crossing, and
thus there are $2^{N}$ states of a diagram with $N$ crossings.
 In  a state we label each smoothing with $A$ or $A^{-1}$ according to the left-right convention 
discussed in Property $3$ (see Figure 14). The label is called a {\em vertex weight} of the state.
There are two evaluations related to a state. The first one is the product of the vertex weights,
denoted  

$$<K|S>.$$
The second evaluation is the number of loops in the state $S$, denoted  $$||S||.$$
  
\noindent Define the {\em state summation}, $<K>$, by the formula 

$$<K> \, = \sum_{S} <K|S>\delta^{||S||-1}.$$
It follows from this definition that $<K>$ satisfies the equations
  
$$<\mbox{\large $\chi$}> \, = A <\mbox{\large $\asymp$}> + A^{-1} <)(>,$$
$$<K \amalg  O> \, = \delta<K>,$$
$$<O> \, =1.$$
  
\noindent The first equation expresses the fact that the entire set of states of a given diagram is
the union, with respect to a given crossing, of those states with an $A$-type smoothing and those
 with an $A^{-1}$-type smoothing at that crossing. The second and the third equation
are clear from the formula defining the state summation. Hence this state summation produces the
bracket polynomial as we have described it at the beginning of the  section. 

\bigbreak

\noindent {\bf Remark.} The bracket polynomial provides a connection between  knot theory and physics, in that the state summation 
expression for it exhibits it as a generalized partition function defined on the knot diagram. Partition functions
are ubiquitous in statistical mechanics, where they express the summation over all states of the physical system of 
probability weighting functions for the individual states. Such physical partition functions contain large amounts of 
information about the corresponding physical system. Some of this information is directly present in the properties of the 
function, such as the location of critical points and phase transition. Some of the information can be obtained by differentiating the 
partition function, or performing other mathematical operations on it. See \cite{KP,Diagrams}.
\bigbreak

\subsection{Thistlethwaite's Example}
View Figure 16. Here we have a version of a link $L$ discovered by Morwen Thistlethwaite \cite{Morwen} in
December 2000. We discuss some theory behind this link in the next subsection. It is a link that is linked but whose linking is not 
detectable by the Jones polynomial. One can verify such properties by using a computer program, or by the algebraic techniques described
below.  

$$ \picill3inby3in(F16) $$
\begin{center} {\bf Figure 16 -- Thistlethwaite's Link} \end{center}
\bigbreak

\subsection{Present Status of Links Not Detectable by the Jones Polynomial}
In this section we give a quick review of the status of our work \cite{EKT}
producing infinite families of distinct links all evaluating as unlinks by the 
Jones polynomial.  
\bigbreak

A tangle (2-tangle) consists in an embedding of two arcs in a three-ball (and possibly some circles
embedded in the interior of the three-ball) such that the endpoints of the arcs are
on the boundary of the three-ball. One usually depicts the arcs as crossing the boundary transversely so that the
tangle is seen as the embedding in the three-ball augmented by four segments emanating from the ball, each from 
the intersection of the arcs with the boundary. These four segments are the {\em exterior edges} of the 
tangle, and are used for operations that form new tangles and new knots and links from given tangles. 
Two tangles in a given three-ball are said to be {\em topologically equivalent} if there is an ambient isotopy
from one to the other in the given three-ball, fixing the intersections of the tangles with the boundary.
\bigbreak

It is customary to illustrate tangles with a diagram that consists in a box  (within which are the arcs of the
tangle) and with the exterior edges emanating from the box  in the NorthWest (NW), NorthEast (NE), SouthWest (SW)
and SouthEast (SE) directions. Given tangles $T$ and $S$, one defines the {\em sum}, denoted $T+S$ by placing
the diagram for $S$ to the right of the diagram for $T$ and attaching the NE edge of $T$ to the NW edge of $S$, and
the SE edge of $T$ to the SW edge of $S$. The resulting tangle $T+S$ has exterior edges corresponding to the 
NW and SW edges of $T$ and the NE and SE edges of $S$. 
There are two ways to create links associated to a tangle $T.$ The {\em numerator} $T^{N}$ is obtained, by attaching
the (top) NW and NE edges of $T$ together and attaching the (bottom) SW and SE edges together. 
The denominator $T^{D}$ is obtained,
by attaching the (left side) NW and SW edges together and attaching the (right side) NE and SE edges together.  
We denote by $[0]$ the tangle with only unknotted arcs (no embedded circles) with one arc connecting, within
the three-ball, the (top points) NW intersection point with the NE intersection point, and the other arc connecting
the (bottom points) SW intersection point with the SE intersection point. 
A ninety degree turn of the tangle $[0]$ produces
the tangle $[\infty]$ with connections between NW and SW and between NE and SE. One then can prove the basic 
formula for any tangle $T$
$$<T> = \alpha_{T} <[0]> + \beta_{T}<[\infty]>$$
\noindent where $\alpha_{T}$ and $\beta_{T}$ are well-defined polynomial invariants (of regular isotopy) of the
tangle $T.$  From this formula one can deduce that 
$$<T^{N}> = \alpha_{T} d + \beta_{T}$$
\noindent and
$$<T^{D}> = \alpha_{T} + \beta_{T} d.$$
\bigbreak

We define the {\em bracket vector} of $T$ to be the ordered pair $(\alpha_{T}, \beta_{T})$ and denote it by
$br(T)$, viewing it as a column vector so that $br(T)^{t} = (\alpha_{T}, \beta_{T})$ where $v^{t}$ denotes the 
transpose of the vector $v.$ With this notation the two formulas above for the evaluation for numerator and 
denominator of a tangle become the single matrix equation

$$\left[
 \begin{array}{c}
      <T^{N}> \\
      <T^{D}>
 \end{array}
 \right] =  \left[
 \begin{array}{cc}
      d & 1 \\
      1 & d
 \end{array}
 \right] br(T).$$
\bigbreak
 
 We then use this formalism to express the bracket polynomial for our examples. The class of examples that 
 we considered are each denoted by $H(T,U)$ where $T$ and $U$ are each tangles and $H(T,U)$ is a satellite
 of the Hopf link that conforms to the pattern shown in Figure 17, formed by clasping together the numerators of the
 tangles $T$ and $U.$ Our method is based on a transformation $H(T,U) \longrightarrow H(T,U)^{\omega}$,
 whereby the tangles $T$ and $U$ are cut out and reglued by certain specific homeomorphisms of the tangle 
 boundaries. This transformation can be specified by a modification described by a specific rational tangle
 and its mirror image. Like mutation, the transformation $\omega$ preserves the bracket polynomial. However,
 it is more effective than mutation in generating examples, as a trivial link can be transformed to a prime
 link, and repeated application yields an infinite sequence of inequivalent links.
\bigbreak

$$ \picill3inby1.5in(F17) $$
\begin{center} {\bf Figure 17 -- Hopf Link Satellite H(T,U)} \end{center}
\bigbreak

$$ \picill5inby4.5in(F18) $$
\begin{center} {\bf Figure 18 -- The Omega Operations} \end{center}
\bigbreak

$$ \picill5inby4in(F19) $$
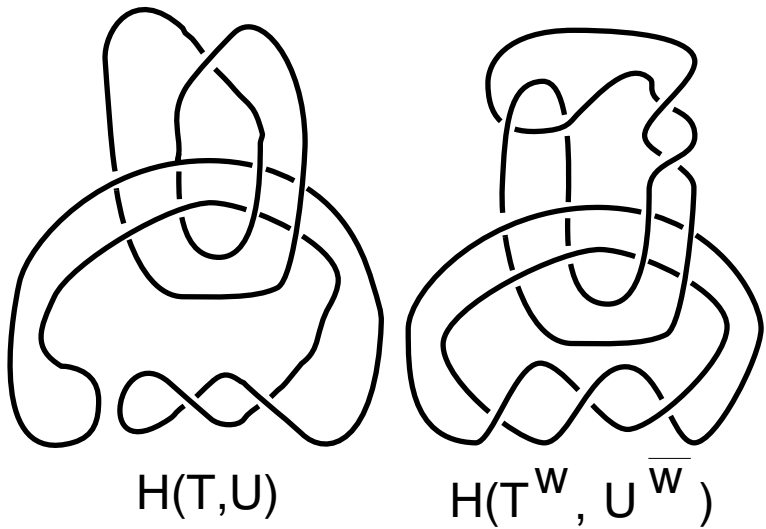
\begin{center} {\bf Figure 19 -- Applying Omega Operations to an Unlink} \end{center}
\bigbreak

 Specifically, the transformation $H(T,U)^{\omega}$ is given by the formula 
 $$H(T,U)^{\omega} = H(T^{\omega}, U^{\bar{\omega}})$$ 
 \noindent where the tangle operations $T^{\omega}$ and  $U^{\bar{\omega}})$ are as shown in Figure 18.
 By direct calculation, there is a matrix $M$ such that 
 $$<H(T,U)> = br(T)^{t} M br(U)$$
 \noindent and there is a matrix $\Omega$ such that 
 $$br(T^{\omega}) = \Omega br(T)$$
 \noindent and 
 $$br(T^{\bar{\omega}}) = \Omega^{-1} br(T).$$
 \noindent One verifies the identity
 $$\Omega^{t} M \Omega^{-1} = M$$
 \noindent from which it follows that $<H(T,U)>^{\omega} = <H(T,U)>.$
 This completes the sketch of our method for obtaining links that whose linking cannot be seen by the Jones polynomial.
Note that the link constructed as $H(T^{\omega}, U^{\bar{\omega}})$ in Figure 19 has the same Jones polynomial as an unlink of
two components. This shows how the first example found by Thistlethwaite fits into our construction.
\bigbreak

\section{From Quantum Topology to Khovanov Homology}

Quantum topology is the study of  invariants of topological objects whose properties 
emulate partition functions in statistical mechanics, quantum mechanics and quantum field theory.
The FKT model of the Alexander-Conway polynomial and the bracket model of the Jones polynomial were the first examples of such
relationships with partition functions. The Jones polynomial was  generalized to a number of other
polynomial invariants of links (the Homflypt and Kauffman polynomials) and to so-called quantum invariants of 
colored links and trivalent framed knotted graphs (the colorings corresponding to the irreducible representations of
the $SU(2)_{q}$ quantum group). These colored invariants were used to create the 
Reshetikhin-Turaev invariants \cite{T,RT} of closed oriented 3-manifolds. Simultaneously, invariants
of three manifolds were discovered by Witten \cite{Witten} in a quantum-field theoretic approach to invariants,
using functional integrals. It was made clear at the heuristic level that the functional integral
approach and the Reshetikhin-Turaev approach were equivalent, but radically different approaches to 
the same structure. The combination of the two approaches led to the notion of topological quantum
field theories (TQFT's) \cite{Atiyah}. A TQFT is a functor from a categories of closed 
surfaces (or higher dimensional manifolds) and their cobordisms, with  appropriate additional
structures, to the category of finite-dimensional vector spaces. This point of view can be used effectively to
capture the properties of the quantum group based invariants and to model the mathematical essence of 
the functional integral approach.
\bigbreak

A new approach to invariants emerged in 1999  when Mikhail Khovanov \cite{Kho} constructed
a homology theory $\cal H_{*}(D)$ defined through diagrams 
$D$ representing  an oriented link $L$ in three-space,
so that the polynomial Euler characteristic  $\chi(\cal{H_{*}(L))}$ of this theory 
is a version of the Jones polynomial of $L,$ and the homology itself carries more topological 
information than the Jones polynomial. Khovanov's construction is an example of  
categorification.
\bigbreak

The Khovanov theory introduces a new concept in quantum topology, namely that it is possible
to transmute state summations to summations of modules and then to 
create topological invariants by taking an appropriate homology or cohomology of these 
modules. These concepts bring the methods of algebraic topology into the subject of quantum 
topology in a fresh way. They suggest that the direct attempt to produce topological partition
functions for four-dimensional topology may have to undergo a significant detour involving 
cohomological methods before it achieves success.
Khovanov homology and its generalizations have proved to be a powerful and useful tool for
low dimensional topology.
\bigbreak   

We end this introduction to Formal Knot Theory with a sketch of the Khovanov homology, as it is directly related to 
the bracket polynomial model. The reader will see that this is indeed a natural outgrowth of the bracket model.  
\bigbreak

Consider the two smoothings of a crossing. It is natural to put then together
as the top and the bottom of a bit of smooth saddle surface as shown in Figure 20. Then there is a topological evolution from one smoothing to 
the other going through the critical level in the saddle where two lines cross one another. The saddle appears as a 
geometrical/topological transformation from one smoothing to the other. 
The question is, how could such a transformation be imaged in an algebraic or combinatorial structure to make deeper versions of the bracket
polynomial invariant? Khovanov's beautiful idea is to regard the transformation from one smoothing to the other as an ingredient for forming the 
differential for a chain complex. 
\bigbreak

$$ \picill5inby6in(F20) $$
\begin{center} {\bf Figure 20 -- Saddle Points and Differentials} \end{center}
\bigbreak

Specifically, one associates an algebra $\cal{A}$ to each loop in a bracket polynomial state $S$ for the knot or link $K.$ Then one takes the tensor 
product $\cal{A}^{\otimes ||S||}$ where $||S||$ is the number of loops in the state $S.$ The $n$-th level of the chain complex $C(K)$ is denoted
$C_{n}(K)$ and is defined to the the direct sum of the modules $\cal{A}^{\otimes ||S||}$ where $S$ is a state with $n$ sites of type $A.$ The differential 
in this complex is designed to take $C_{n+1}(K)$ to  $C_{n}(K).$ The differential is built from the two possibilities for resmoothing a site, as
illustrated  in Figure 20. In one case two curves become a single curve. In the other case a single curve is transformed into two curves. We have denoted
these two possibilities by $m$ and $\Delta$ respectively. Since each curve is associated to a copy of the algebra $\cal{A}$, it follows that there should be 
maps of algebras $m: \cal{A} \otimes \cal{A} \longrightarrow \cal{A}$ and $\Delta: \cal{A} \longrightarrow \cal{A} \otimes \cal{A}.$ Thus the algebra 
needs a product $m$ and a coproduct $\Delta$. The differential always goes from a smoothing of type $A$ to a smoothing of type $A^{-1}$ when restricted 
to a given loop, or pair of loops in a state. Since a resmoothing of a type $A$ state changes it to a type $A^{-1}$ state, the number of type $A$ states
is reduced by one when the differential is performed. The total differential from  $C_{n+1}(K)$ to  $C_{n}(K)$ is assembled from a signed summation 
over all single-site maps as described above. The signs require ordering choices and we refer the reader to \cite{Dror} for more information (where the
complex is expressed in cohomology rather than homology). 
\bigbreak

A key point about this construction is illustrated in Figure 21. In order that the differential in the complex be well-defined and have the property
that $d^{2} = 0,$ one needs that different orders of application of it will coincide. It is easiest to see this necessary compatibility condition by 
thinking of the elementary bit of differential as $m$ or $\Delta$, depicted as a bit of saddle surface beginning with two loops and ending with one loop, or
beginning with one loop and ending with two loops. Then the compatibility corresponds to the equality of the compositions $F=G=H$ indicated in Figure 21 as 
$$F = \Delta \circ m,$$
$$G = (m \otimes 1)\circ(1 \otimes \Delta),$$
$$H = (1\otimes \Delta)\circ(m \otimes 1).$$
These are the conditions for $\cal{A}$ to be a Frobenius algebra \cite{Dror}.
For an example of an algebra satisfying these requirements, take $\cal{A}$ to have basis $\{ 1, X \}$ with $1$ acting as a multiplicative identity and 
$$X^{2} = m(X \otimes X) = 0,$$ 
$$\Delta(1) = 1 \otimes X + X \otimes 1,$$
and $$\Delta(X) = X \otimes X.$$ This algebra satisfies the requirements and produces 
an homology theory that is invariant under the Reidemeister moves (when one takes into account grading changes that the moves can produce).
\bigbreak

$$ \picill3inby4in(F21) $$
\begin{center} {\bf Figure 21 -- The Frobenius Algebra Conditions} \end{center}
\bigbreak
 
Khovanov homology grows naturally out of the bracket state sum,creating a new technique for obtaining topological
information from the states of the link diagram. 
\bigbreak 

\subsection {Knot Floer Homology}
Finally, we mention the work of Ozsv\'{a}th and Szab\'{o} \cite{OS} on Heegard Floer homology. This is a homology theory associated to knots that categorifies
the  Alexander Polynomial. The chain complex for this theory has a basis that is in one-to-one correspondence with the states in the Formal Knot Theory 
model for the Alexander-Conway polynomial. At this writing, there is no strictly combinatorial description for the differentials in the chain complex for this
homology theory. The theory appears to depend upon more complex geometry than the Khovanov invariant. 
\bigbreak


\begin{thebibliography}{99}

\bibitem{Alex} J. W. Alexander, Topological invariants of knots and links, {\it Trans. Amer. Math. Soc.}, Vol. 30, Issue 2, (April 1928), 275-306.

\bibitem{Atiyah} M.F. Atiyah, {\em The Geometry and Physics of Knots}, Cambridge
University Press, 1990.


\bibitem{Dror}
 D. Bar-Natan, On Khovanov's categorification of the Jones polynomial. {\em
Algebraic and Geometric Topology}, Vol. 2 (2002), pp. 337-370.


\bibitem{Conway} J. H. Conway, An enumeration of knots and links and some of
their algebraic properties, in {\em Computational Problems in Abstract Algebra},
Pergammon Press, N.Y.,1970, pp. 329-358.


\bibitem{Fox1} R. H. Fox, A quick trip through knot theory. In ``Topology of $3$-Manifolds",
ed. by M. K. Fort Jr., Prentice Hall (1962), 120-167.

\bibitem{Fox2} R. H. Crowell and R. H. Fox, ``Introduction to Knot Theory", Blaisdell Pub. Co. (1963).

\bibitem{GL} P. M. Gilmer and R. A. Litherland, The duality conjecture in formal knot theory, {\it Osaka J. Math.} Vol. 23, (1986), 275-306.


\bibitem{Jones} V.F.R. Jones, A polynomial invariant of links via von Neumann
algebras, {\em Bull. Amer. Math. Soc.}, 1985, No. 129, pp. 103-112.


\bibitem{Kauffman-Conway} L.H. Kauffman, The Conway polynomial, {\em Topology},
{\bf 20} (1980), pp. 101-108.

\bibitem{FKT} L. H. Kauffman, {\em Formal Knot Theory}, Princeton
University Press, Lecture Notes Series  30 (1983).

\bibitem{OnKnots} L.H. Kauffman, {\em On Knots}, Annals Study No. 115,
{\em Princeton University Press} (1987)

\bibitem{Bracket} L.H.Kauffman, State Models and the Jones Polynomial, {\em
Topology},Vol. 26, 1987,pp. 395-407.

\bibitem{New} L.H.Kauffman, New invariants in the theory of knots,
{\em Amer. Math. Monthly}, Vol.95,No.3,March 1988. pp 195-242.

\bibitem{KP}
 L. H. Kauffman, ``Knots and Physics", World Scientific,
Singapore/New Jersey/London/Hong Kong, 1991, 1994, 2001.


\bibitem{TL} L. H. Kauffman and S. L. Lins, {\em Temperley-Lieb Recoupling Theory
and Invariants of 3- Manifolds}, Annals of Mathematics Study 114, Princeton Univ.
Press,1994.

\bibitem{Diagrams} L. H. Kauffman, math.GN/0410329, Knot diagrammatics. "Handbook of Knot Theory``, edited by Menasco and
Thistlethwaite,  233--318, Elsevier B. V., Amsterdam, 2005.

\bibitem{Kho}
 Mikhail Khovanov, A categorification of the Jones polynomial, {\it Duke Math. J.} 101 (2000), no. 3, 359--426.


\bibitem{Rolfsen} D. Rolfsen, ``Knots and Links", Publish or Perish Press (1987).

\bibitem{Reidemeister} K. Reidemeister, {\em Knotentheorie}, Chelsea Pub. Co.,
New York, 1948, Copyright 1932, Julius Springer, Berlin.

\bibitem{T}
M.B. Thistlethwaite,
A spanning tree expansion of the Jones polynomial,
{\em Topology} {\bf 26} no. 3 (1987), 297--309.

\bibitem{Morwen}
M. Thistlethwaite, Links with trivial Jones polynomial, JKTR, Vol. 10, No. 4 (2001), 641-643.

\bibitem{EKT} S. Eliahou, L. Kauffman and M. Thistlethwaite, Infinite families of links with trivial Jones polynomial,
{\em Topology}, {\bf 42}, pp. 155--169.

\bibitem{OS} P. Ozsv\'{a}th and Z. Szab\'{o}, Heegard Floer homology and alternating knots, {\em Geometry and Topology}, Vol. 7 (2003), 225-254.

\bibitem{Propp} J. Propp, Lattice structure for orientations of graphs, MIT (1993), (unpublished).

\bibitem{RT} N.Y. Reshetikhin and V. Turaev, Invariants of Three-Manifolds via
link polynomials and quantum groups, {\em Invent. Math.},Vol.103,1991, pp.
547-597.

\bibitem{T} V.G.Turaev, The Yang-Baxter equations and invariants of links,
LOMI preprint E-3-87, Steklov Institute, Leningrad, USSR., {\em Inventiones
Math.},Vol. 92, Fasc.3, pp. 527-553.


\bibitem{Witten} Edward Witten, Quantum field Theory and the Jones Polynomial,
{\em Commun. Math. Phys.},vol. 121, 1989, pp. 351-399.




\end{thebibliography}
\end{document}